\nopagenumbers

\headline={\hss\tenrm\folio\hss}
\magnification
\magstep1
\baselineskip 18pt
\hsize=175truemm
\vsize=240truemm

\parindent 0.0cm
\centerline {\bf KOEBE DOMAIN OF STARLIKE FUNCTIONS OF COMPLEX }
\centerline {\bf ORDER WITH MONTEL NORMALIZATION }
\vskip 0.8cm
\centerline {  YA\c SAR POLATO\u GLU $^{*}$,MET\.IN BOLCAL $^{**}$ AND ARZU \c SEN $^{***}$}
\vskip 0.5cm
\par {\bf ABSTRACT:}    Let $\> S^{*}(1-b)\> $($b \not= 0 $ complex) denote the class of functions
 $ f(z)=z+\alpha_{2}z^{2}+\cdots  $ analytic in $ D={z \mid \vert z \vert < 1 } $
which satisfies,for $z=e^{i\theta } \in D $,$( f(z)/z )\not=  0 $ in $D$,and
$$ Re \Biggr [ 1+ {1 \over b} \Biggr ( z { { f^{'}(z)} \over f(z) }-1 \Biggl ) \Biggl ] > 0 .$$
The aim of this paper is to give the Koebe domain of the above mentioned class.
\vskip 0.5cm
\par {\bf INTRODUCTION:} Let F denote the class of function
$\> f(z)=z+\alpha_{2}z^{2}+\cdots \> $ which are analytic in $D$ .We shall need the following definitions.
\vskip 0.3cm
\par {\bf DEFINITIONS.1.1:}A function $f(z) \in $ F is said to be starlike function
of complex order $(1-b)$,($b \not= 0 $ complex) ,that is $f(z) \in S^{*}(1-b)$ if and only if
$( f(z)/z )\not=  0 $ in $D$,and
$$   Re \Biggr [ 1+ {1 \over b} \Biggr ( z { { f^{'}(z)} \over f(z) }-1 \Biggl ) \Biggl ] > 0 \quad,\quad z\in D  \eqno(1.1) $$
It should be noticed that by giving specific values to $b$ we obtain the following
important subclasses.[5]

\par I)$\>\>\>\> b=1,S^{*}(1-b)=S(0)=S^{*}$ is the well known classes of starlike functions.
\par II)$\>\>\> b=1-\alpha , (0 \leq \alpha<1),S^{*}(1-b)=S(\alpha)$ is the class of starlike functions of order $\alpha $.
\par III)$\>b= \cos \lambda .e^{-i\lambda},S^{*}(1-b)=S^{*}(1-\cos \lambda .e^{-i\lambda})=S_{\lambda }^{*},
\vert \lambda \vert < {{\pi } \over 2 }$ is the class of spirallike functions.
\par IV)$\>b=1-\alpha \cos \lambda .e^{-i\lambda},\vert \lambda \vert < {{\pi } \over 2 },
0 \leq \alpha<1,
S^{*}(1-b)=S^{*}(1-(1-\alpha )\cos \lambda .e^{-i\lambda})=S_{\lambda,\alpha }^{*}$ is the class of spirallike functions of order $\alpha $.
\vskip 0.3cm
\par {\bf DEFINITIONS 1.2:}The Koebe domain for the family $F$ is denoted by $K(F)$ and by the definition this is largest domain is
contained in $f(D)$ for every function $f(z)$ in $F$.From this definition the Koebe domain has the following properties.
\par I) $K(F)$ is the collection of points $w$ such that $w$ is in $f(D)$for every function $f(z)$ in $F$,symbols
$$K(F) = \bigcap\limits_{f(z)\in F}f(D) \eqno (1.2)$$
\par II)Supposing the set $F$ is invariant under the rotation,so $e^{i\alpha}.f(e^{-i\alpha}.z)$
is whenever $f(z)$ in $F$.
\vskip 0.3cm
\par Then the Koebe Domain will be either the single point $w=0$ or an open disc.$\vert w \vert <R$.In the second case $R$ is often easy to find indeed,
supposing that we have a sharp lover bound $M(r)$ for $f(re^{i\theta})$for all functions
in $F$,and $F$ contains only
$$ R=\lim\limits_{r\rightarrow 1^{-}}M(r) \eqno(1.3) $$,
\par gives the disc $\vert w \vert <R$ as the Koebe Domain for the set $F$.[see1.2]
\vskip 0.3cm
\par  MONTEL TYPE NORMALIZATION: We can also impose a Montel type normalization.
This means that  for some fixed $r_{0}$with $0<r_{0}<1$,we consider the family of functions $f(z)$ regular and univalent in $D$
with $f(0)=0,f(0)=1,f(r_{0})=r_{0}$.
\vskip 0.3cm
We note that: If the class of starlike functions is normalized by the Montel type normalization,then the
class is denoted by.
$$ S^{*}_{montel}(1-b)$$
\par {\bf II.KOEBE DOMAIN FOR THE CLASS OF STARLIKE FUNCTIONS OF COMPLEX ORDER WITH MONTEL NORMALIZATION.}
\par In this section we shall give the Koebe Domain for the class of starlike functions of complex order with Montel normalization.
\par {\bf THEOREM:}Let $f(z) \in S^{*}(1-b)$,then
$$\eqalign {
{
{2\vert u \vert \Big (1-\vert v \vert ^{2}\Big )^{2b}} \over
{(1+\vert b \vert)\vert v \vert \vert 1-u\overline v \vert^{2b-2}
\Big[ \vert 1-u\overline v \vert + \vert u-v \vert \Big ] ^{2} }} &\leq 
 \Big \vert {{f(u)} \over {f(v)}} \Big \vert \leq \cr \cr
&\leq {{2\vert u \vert \Big (1-\vert v \vert ^{2}\Big )^{2b}} \over
{(1+\vert b \vert)\vert v \vert \vert 1-u\overline v \vert^{2b-2}
\Big[ \vert 1-u\overline v \vert - \vert u-v \vert \Big ] ^{2} }}  }$$
\par holds for $u,v \in D,u\not = v$.This bound is sharp,because the extremal function is $w=f(z)$ defined for $\vert z \vert <1$ by
$$ w= f_{*}(z)={z \over {(1-z)^{2b}}}$$
\par  Proof: We consider the M"bius transformation
$$ u= {{z+v} \over {1+\overline v z}}\Leftrightarrow z={{u-v} \over {1-u \overline v}} \eqno (2.1)$$
\par which is analytic and univalent in $D$ ,and this M"bius transformation maps the unit
disc on to the itself.Now we define the function:
$$F(z)={ {\alpha z f \Bigl ( {{z+\alpha }\over {1+\overline \alpha z}} \Bigr ) } \over
{f(z)(z+\alpha)(1+\overline \alpha z)^{2b-1}}}\Leftrightarrow F(z)={
{v(1-u\overline v)^{2b-1}(u-v)f(u)}\over{u(1-\vert v \vert^{2})^{2b}f(v)}},\alpha=v \eqno(2.2)$$
\par This functions starlike function of complex order in the unit disc$D$.On the other hand if $g(z)$ is starlike fnction of complex order,then
$${ {2\vert z \vert } \over {(1+\vert b \vert ).(1+\vert z \vert )^{2}} }\leq
F(z) \leq 
{ {2\vert z \vert } \over {(1+\vert b \vert ).(1-\vert z \vert )^{2}} } \eqno (2.3)$$
\par holds[5].Therefore applying the inequality(2.3) to the function $F(z)$ which is defined by (2.2) and using the relation (2.1) we get
$$
{{2 \Bigl \vert {{u-v} \over {1-u\overline v}}\Bigr \vert} \over
{(1+\vert b \vert ). \Bigl ( 1+\Bigl \vert 
{{u-v} \over {1-u\overline v}}\Bigr \vert \Bigr)^{2}}}\leq
\Biggl \vert {{v(1-u\overline v)^{2b-1}(u-v)f(u)} \over {u(1-\vert v \vert ^{2})^{2b}f(v)}}\Biggr \vert\leq
{{2 \Bigl \vert {{u-v} \over {1-u\overline v}}\Bigr \vert} \over
{(1+\vert b \vert ). \Bigl ( 1-\Bigl \vert 
{{u-v} \over {1-u\overline v}}\Bigr \vert \Bigr)^{2}}} \eqno (2.4)$$
\par Simple calculations from (2.4) show that this theorem is true.
\vskip 0.3cm
\par COROLLARY.2.1.The Koebe Domain of starlike function of complex order with Montel
type nornalization is:
$$R={
{(1-r_{0}^{2})^{2b}} \over {2(1+\vert b \vert )(1-2r_{0}\cos\theta +r_{0}^{2})^{2b}}},0\leq \theta \leq 2\pi $$
\par Indeed;If we take $v=f(r_{0})=r_{0},0<r_{0}<1,u=z=re^{i\theta }$ in theorem (2.1) we obtain.
$$\eqalignno {& 
{
{2 \vert re^{i\theta } \vert (1-r_{0}^{2})^{2b} }
\over
{(1+\vert b \vert ). {\vert 1-r_{0}re^{i\theta } \vert }^{2b-2}
\Bigl [ \vert 1-r_{0}re^{i\theta } \vert + \vert re^{i\theta }-r_{0} \vert  \Bigr ]^{2} }
}\leq      \cr \cr
&\leq  \vert f(z) \vert \leq \cr \cr
&\leq {
{2\vert re^{i\theta } \vert (1-r_{0}^{2})^{2b} }
\over
{(1+\vert b \vert ). {\vert 1-r_{0}re^{i\theta } \vert }^{2b-2}
\Bigl [ \vert 1-r_{0}re^{i\theta } \vert - \vert re^{i\theta }-r_{0} \vert  \Bigr ]^{2} } } &(2.5) \cr } $$
\par Therefore we have
$$M(r) =
{
{2 \vert re^{i\theta } \vert (1-r_{0}^{2})^{2b} }
\over
{(1+\vert b \vert ). {\vert 1-r_{0}re^{i\theta } \vert }^{2b-2}
\Bigl [ \vert 1-r_{0}re^{i\theta } \vert + \vert re^{i\theta }-r_{0} \vert  \Bigr ]^{2} }} $$
\par if we take $\lim\limits _{r\rightarrow 1^{-}}M(r)$ we get
$$R=K(S^{*}_{montel}(1-b))={
{(1-r_{0}^{2})^{2b}} \over
{2(1+\vert b \vert ).(1-2r_{0}\cos \theta + r_{0}^{2})^{b}}}\eqno (2.6)  $$
\par Giving the specific values to $b$ we obtain the following results.
\par (I) For $b=1$

$$ R=K \Bigl (S^{*}_{montel } \Bigr )={
{(1-r_{0}^{2})^{2} } \over
{4(1-2r_{0}\cos \theta + r_{0}^{2}) }} $$

\par (II) For $b=1-\alpha \>\>,\>\> 0 \leq \alpha <1 $

$$ R=K \Bigl (S^{*}_{montel } (\alpha )\Bigr )={
{(1-r_{0}^{2})^{2(1-\alpha )} } \over
{2\alpha(1-2r_{0}\cos \theta + r_{0}^{2})^{1-\alpha }}} $$

\par (III) For $b=\cos\lambda e^{-i\lambda } \>\>,\>\>(\vert \lambda \vert < \pi / 2)$
$$ R=K \Biggl (S^{*}_{montel } \Big ( 1-\cos\lambda .e^{-i \lambda } \Big ) \Biggr )
={
{(1-r_{0}^{2})^{2\cos\lambda . e^{-i \lambda }} } \over
{2(1+\cos\lambda } )(1-2r_{0}\cos \theta + r_{0}^{2})^{\cos\lambda .e^{-i \lambda} }} $$
\par (IV) For $b=(1-\alpha )
\cos \lambda .e^{-i\lambda } \>\>,\>\>(\vert \lambda \vert < \pi / 2)\>\>,\>\> 0 \leq \alpha < 1 $
$$ R=K \Biggl (S^{*}_{montel } \Big ( 1-(1-\alpha )\cos\lambda .e^{-i \lambda } \Big ) \Biggr )
={
{(1-r_{0}^{2})^{(1-\alpha )\cos\lambda .e^{-i \lambda }} } \over
{(1+(1-\alpha )\cos\lambda } )(1-2r_{0}\cos \theta + r_{0}^{2})^{(1-\alpha )\cos\lambda .e^{-i \lambda} }} $$
\vskip 0.8cm
\par \centerline { R E F E R E N C E S }
\par [1.]P.L.Duren,Univalent functions,Springer Verlag,New York 1981
\par [2.]A.W.Goodman,Univalent functions,Volume I and Volume II,Mariner Publ.Co.Inc.Tampa
\par Florida 1983
\par [3.]J.Krzyz,On a problem of P.Montel,Ann.Polon.Math.,12.(1962) 55-60
\par [4.]Ky Fan,Distortion of Univalent functions,J.Math.Anal.Appl.66(1979)p:626-631
\par [5.]Polato§lu Y,Koebe Domain for convex functions of complex order,The Punjab University Journal of Mathematics.Vol.XXX(1979)(To appear)
\vskip 0.5cm
\par [$^{*}$.]YAžAR POLATO¦LU:Department of Mathematics,Faculty of Science,University of Kltr.
\par [$^{**}$.]MET˜N BOLCAL:Department of Mathematics,Faculty of Science,University of Kltr.
\par [$^{***}$.]ARZU žEN:Department of Mathematics,Faculty of Science,University of Kltr.

\end